\documentclass[A4paper,12pt]{book}
\usepackage{amssymb,amsmath}
\usepackage[russian]{babel}
\usepackage{graphicx}
\newcount\exnom\exnom=0

\long\def\exo/{\vspace{0.1cm} \noindent\advance\exnom by1{
{\the\exnom}.}}
\newcommand{\z}{\exo/  }
\newcommand{\ds}{\displaystyle}

\begin{document}
\thispagestyle{empty}

\begin{center}  \textbf{ PROJECTIONS OF  $\mathbf{k}$ -- SIMPLEX ONTO THE SUBSIMPLICES OF ARBITRARY TYPE ARE  DERIVATIONS}
\end{center}

\begin{center} \textbf{Dimitrinka Vladeva}
\end{center}

\vspace{2mm}

\noindent\emph{\textbf{Abstract:}
The aim of this paper is to prove that there is a  projection of  an arbitrary $k$ --simplex onto  $(m-\ell)$ -- subsimplex, where $ 1 \leq \ell < m \leq k -1$, which is a derivation.}

\vspace{1mm}


\vspace{5mm}

{\bf 1 \hspace{1mm} Introduction and preliminaries}

\vspace{3mm}

Differential algebra is an area of algebra in the study of algebraic structures equipped with one or finitely many
derivations that are additive maps satisfied Leibniz product rule.
A simple examples are  the usual derivatives on various algebras consisting of differentiable functions in sense of analysis.
In 1950, Ritt [5] and in 1973, Kolchin [3] wrote their  classical  books on differential algebra.
During the last few decades there has been a great deal of works concerning derivations in fields and rings, in Lie rings, in skew polynomial rings and other algebraic structures.

About derivations in semirings we know:

--- the definition in Golan's book [1],

--- properties of derivations in strings and some simplicial complexes of strings -- [9],

---  the projections on the strings or on the triangles of 2, 3 and 4 simplices which are derivations --  [11] -- [14],

---  the projections on the triangles $\triangle^{(n)}\{a_0,a_1,a_m\}$, where
$2 \leq m \leq k-1$, of an arbitrary simplex $\sigma^{(n)}\{a_0,a_1,\ldots,a_{k-1}\}$ which are derivations -- [15],

---  the projections on the all basic strings of  an arbitrary simplex $\sigma^{(n)}\{a_0,a_1,\ldots,a_{k-1}\}$ which are derivations -- [16].

The theorem proved here generalise most of these results.

\vspace{3mm}

The endomorphism semirings of a finite semilattice are well-established, see
[2], [7], [8] and [17]. Concerning the background of  simplicial complexes, algebraic topology and combinatorics a reader can  refer to [4], and [6].

\vspace{3mm}

 An algebra $R = (R,+,.)$ with two binary operations $+$ and $\cdot$ on $R$, is called a {\emph{semiring}} if:

\textbf{1.} $\;(R,+)$ is a commutative semigroup, \hspace{2mm}

\textbf{2.} $\; (R,\cdot)$ is a semigroup, \hspace{4mm}

\textbf{3.}   distributive laws hold $ x\cdot(y + z) = x\cdot y + x\cdot z$ and $(x + y)\cdot z = x\cdot z + y\cdot z$
for any $x, y, z \in R$.

For a join-semilattice $\left(\mathcal{M},\vee\right)$  set $\mathcal{E}_\mathcal{M}$ of the endomorphisms of $\mathcal{M}$ is a semiring
 with respect to the addition and multiplication defined by:

 $\bullet \; h = f + g \; \mbox{when} \; h(x) = f(x)\vee g(x) \; \mbox{for all} \; x \in \mathcal{M}$,

 $\bullet \; h = f\cdot g \; \mbox{when} \; h(x) = f\left(g(x)\right) \; \mbox{for all} \; x \in \mathcal{M}$.

 This semiring is called the \emph{endomorphism semiring} of $\mathcal{M}$.

\vspace{2mm}

In this article, all semilattices are finite chains. Following [7] we fix a finite chain $\mathcal{C}_n = \; \left(\{0, 1, \ldots, n - 1\}\,,\,\vee\right)\;$
and denote the endomorphism semiring of this chain with $\widehat{\mathcal{E}}_{\mathcal{C}_n}$. We do not assume that $\alpha(0) = 0$ for arbitrary
$\alpha \in \widehat{\mathcal{E}}_{\mathcal{C}_n}$. So, there is not a zero in  endomorphism semiring $\widehat{\mathcal{E}}_{\mathcal{C}_n}$.

\vspace{2mm}

 In [10] there is a new treatment of the subsemirings of  endomorphism
 semiring $\widehat{\mathcal{E}}_{\mathcal{C}_n}$ of a finite chain.
 For arbitrary elements $a_0, a_1, \ldots, a_{k-1} \in \mathrm{C}_n$, where $k \leq n$ and $a_0 < a_1 < \ldots < a_{k-1}$ we denote $A = \{a_0, a_1, \ldots, a_{k-1}\}$.
Now, consider  endomorphisms $\alpha \in \widehat{\mathcal{E}}_{\mathcal{C}_n}$ with $Im(\alpha) \subseteq A$.  The set of the all such endomorphisms $\alpha$  is a maximal simplex. We denote this simplex  by
 $\sigma^{(n)}_k(A) = \sigma^{(n)}\{a_0, a_1, \ldots, a_{k-1}\}$.

The  endomorphisms $\alpha \in \sigma^{(n)}\{a_0, a_1, \ldots, a_{k-1}\}$ such that
$$\alpha(0) = \cdots = \alpha(i_0-1) = a_0, \alpha(i_0) = \cdots = \alpha(i_0+ i_1 -1) = a_1, \; \cdots$$ $$\alpha(i_0 + \cdots +i_{k-2}) = \cdots = \alpha(i_0 + \cdots + i_{k-1} - 1) = a_{k-1}$$
 we denote  by $\alpha = (a_0)_{i_0}(a_1)_{i_1} \ldots (a_{k-1})_{i_{k-1}}$, where
 $\ds \sum_{p=0}^{k-1} i_p = n$.

\vspace{5mm}

{\bf 2 \hspace{1mm} Addition and multiplication of endomorphisms}
\vspace{3mm}

Here we show how to signify the sum of two endomorphisms.

\vspace{5mm}

\textbf{Lemma 1.} \textsl{Let $\alpha, \beta \in \sigma^{(n)}\{a_0,\ldots,a_{k-1}\}$, where
 $\alpha = (a_0)_{i_0} \cdots (a_{k-1})_{i_{k-1}}$, $\beta = (a_0)_{j_0} \cdots (a_{k-1})_{j_{k-1}}$ and $\ds \sum_{p=0}^{k-1} i_p = \sum_{p=0}^{k-1} j_p = n$. Let  $\gamma = (a_0)_{h_0} \cdots (a_{k-1})_{h_{k-1}}$, where}
$$h_0 = \min\{i_0,j_0\},$$
$$h_s = \min\left\{ \sum_{p=0}^{s} i_p -  \sum_{p=0}^{s-1} h_p\, , \, \sum_{p=0}^{s} j_p -  \sum_{p=0}^{s-1} h_p\right\}, \forall s = 1, \ldots k-1. \eqno{(1)}$$

Then $\gamma = \alpha + \beta$.

\vspace{2mm}

\emph{Proof.} From (1) follows
$$\sum_{p=0}^{s} h_p  = \min\left\{ \sum_{p=0}^{s} i_p \, , \, \sum_{p=0}^{s} j_p \right\}, \forall s = 1, \ldots k-1. $$

Now for $s = k-1$ we find
$$\sum_{p=0}^{k-1} h_p  = \min\left\{ \sum_{p=0}^{k-1} i_p \, , \, \sum_{p=0}^{k-1} j_p \right\} = \min\{n,n\} = n. $$
Hence, $\gamma \in \sigma^{(n)}\{a_0,\ldots,a_{k-1}\}$.

Let $t = 0, 1, \ldots, k - 1$. Let for some $s$ have $\ds \sum_{p = 0}^{s} i_p < t \leq \sum_{p = 0}^{s+1} i_p$. Then it follows $\alpha(t) = a_{s+1}$.

Let us suppose that $\ds t \leq \sum_{p = 0}^{s} j_p$. Then it follows $\ds \sum_{p = 0}^{s} h_p = \sum_{p = 0}^{s} i_p$. Since
$\ds  t \leq \sum_{p = 0}^{s+1} i_p$ and $\ds t \leq \sum_{p = 0}^{s} j_p \leq \sum_{p = 0}^{s+1} j_p$, it follows $\ds  t \leq \sum_{p = 0}^{s+1} h_p$. So $\ds \sum_{p = 0}^{s} h_p < t \leq \sum_{p = 0}^{s+1} h_p$, that is $\gamma(t) = a_{s+1} = \alpha(t)$.

On the other hand, from $\ds  t \leq \sum_{p = 0}^{s} j_p$ it follows that $\beta(t) \leq a_{s}$.

 Thus $\gamma(t) = \max\{\alpha(t),\beta(t)\}$.

 Let us suppose that $\ds \sum_{p = 0}^{s} j_p < t \leq \sum_{p = 0}^{s+ 1} j_p$. Then $\ds \sum_{p = 0}^{s} h_p < t \leq \sum_{p = 0}^{s+ 1} h_p$ and now $\gamma(t) = \beta(t) = \alpha(t)$.

 Let us suppose that $\ds \sum_{p= 0}^{s+1} j_p < t$. Then $\ds \sum_{p= 0}^{s+1} j_p < \sum_{p= 0}^{s+1} i_p$, so, $\ds \sum_{p= 0}^{s+1} h_p = \sum_{p= 0}^{s+1} j_p$. Now $\gamma(t) = \beta(t) \geq a_{s+1}$ and also $\gamma(t) = \max\{\alpha(t),\beta(t)\}$.

 Hence, we prove that $\gamma(t) = \alpha(t) + \beta(t)$.

\vspace{4mm}

For multiplication of two endomorphisms  $\alpha = (a_0)_{i_0} \cdots (a_{k-1})_{i_{k-1}}$ and $\beta = (a_0)_{j_0} \cdots (a_{k-1})_{j_{k-1}}$, where $\ds \sum_{p=0}^{k-1} i_p = \sum_{p=0}^{k-1} j_p = n$,  we consider the following notations.

Let $\beta(a_0) = a'_0$ and $n_0$ be the greatest number such that $$\beta(a_0) = \cdots = \beta(a_{n_0}) = a'_0.$$
Let $a_{n_0+1}$ be the next element in $\{a_0, \ldots, a_{k-1}\}$ after $a_{n_0}$ and $\beta(a_{n_0+1}) = a'_1$. Let $n_1$ be the greatest  number such that
$$\beta(a_{n_0+1}) = \cdots = \beta(a_{n_1}) = a'_1.$$

So, we define a partition of the set $\{a_0, \ldots, a_{k-1}\}$ with properties
$$\beta(a_p) = a'_0, \forall\, p,\; 0 \leq p \leq n_0,$$
$$\beta(a_p) = a'_1, \forall\, p,\; n_0 + 1 \leq p \leq n_1,$$
$$ ............................................................$$
$$\beta(a_p) = a'_r, \forall\, p,\; n_{r-1} + 1 \leq p \leq n_r = k - 1.$$

\vspace{2mm}

Since $(\alpha\beta)(t) = \beta(\alpha(t))$ for $t = 0, \ldots, k - 1$, it follows
$$\alpha\beta = (a'_0)_{\sum_{p=0}^{n_0} i_p}(a'_1)_{\sum_{p=n_0+1}^{n_1} i_p} \cdots (a'_r)_{\sum_{p=n_{r-1}+1}^{n_r} i_p}. \eqno{(2)}$$

Let $\beta(a_{k-1}) = a'_{k-1}$. Since $n_r = k-1$, it follows that $a'_{k-1} = a'_r$. So, $n_{r-1} + 1$  is the smallest number such that
$$\beta(a_{k-1}) = \cdots = \beta(a_{n_{r-1}+1}) = a'_{k-1}.$$

 Then, from (2), it follows
$$\alpha\beta = (a'_{0})_{\sum_{p=0}^{n_0} i_p} \cdots (a'_{u})_{\sum_{p=n_{u-1} + 1}^{n_u} i_p} \cdots  (a'_{k-1})_{\sum_{p=n_{r-1}+1}^{k-1} i_p}. \eqno{(3)}$$
In (3) we note that $1 \leq u \leq k-1$, $\;n_{u-1} + 1$ is the smallest number and ${n_u}$ is the greatest number such that
$\beta(a_{n_{u-1} + 1}) = \cdots = \beta(a_{n_u}) = a'_u$.

\vspace{5mm}

{\bf 3 \hspace{1mm} Projection onto   $(m - \ell)$ -- simplex}

\vspace{3mm}

Let us consider the map
$$\partial_{{m-\ell}}^{k-1 - 0} = \partial_{{m-\ell}}^{k-1} : \sigma^{(n)}\{a_0,\ldots,a_{k-1}\} \rightarrow \sigma^{(n)}\{a_{\ell},\ldots,a_{m}\},$$
where $0 \leq \ell < m \leq k-1$ and $a_\ell, \ldots, a_m$ are consecutive elements of the set $\{a_0,\ldots,a_{k-1}\}$,
so that for any endomorphism $\alpha \in  \sigma^{(n)}\{a_0,\ldots,a_{k-1}\}$,   $\alpha = (a_0)_{i_0} \cdots (a_{k-1})_{i_{k-1}}$, where $\ds \sum_{p=0}^{k-1} i_p = n$,
 $$\partial_{{m-\ell}}^{k-1}(\alpha) = (a_{\ell})_{\sum_{p=0}^{\ell} i_{p}}(a_{\ell +1})_{i_{\ell+1}} \cdots (a_{m-1})_{i_{m-1}}
 (a_m)_{\sum_{p=m}^{k-1} i_{p}}.$$

\vspace{4mm}

\textbf{Lemma 2.} \textsl{For any endomorphisms $\alpha, \beta \in \sigma^{(n)}\{a_0,\ldots,a_{k-1}\}$, it follows}
 $$\partial_{{m-\ell}}^{k-1}(\alpha + \beta) = \partial_{{m-\ell}}^{k-1}(\alpha) + \partial_{{m-\ell}}^{k-1}(\beta), $$
where $0 \leq \ell < m \leq k-1$.

\vspace{2mm}

\emph{Proof.} Let $ \alpha = (a_0)_{i_0} \cdots (a_{k-1})_{i_{k-1}}$ and $\beta = (a_0)_{j_0} \cdots (a_{k-1})_{j_{k-1}}$. On the basis of Lemma~1 we can conclude that
$\alpha + \beta = (a_0)_{h_0} \cdots (a_{k-1})_{h_{k-1}}$, where
$$h_0 = \min\{i_0,j_0\},$$
$$h_s = \min\left\{\sum_{p=0}^{s} i_p - \sum_{p=0}^{s-1} h_p\,,\,\sum_{p=0}^{s} j_p - \sum_{p=0}^{s-1} h_p\right\}\; \mbox{for}\; s = 1, \ldots k-1.$$

Now we consider
 $$\partial_{{m-\ell}}^{k-1}(\alpha) = (a_{\ell})_{\sum_{p=0}^{\ell} i_{p}}(a_{\ell +1})_{i_{\ell+1}} \cdots (a_{m-1})_{i_{m-1}}
 (a_m)_{\sum_{p=m}^{k-1} i_{p}},$$
 $$\partial_{{m-\ell}}^{k-1}(\beta) = (a_{\ell})_{\sum_{p=0}^{\ell} j_{p}}(a_{\ell +1})_{j_{\ell+1}} \cdots (a_{m-1})_{j_{m-1}}
 (a_m)_{\sum_{p=m}^{k-1} j_{p}},$$
 $$\partial_{{m-\ell}}^{k-1}(\alpha + \beta) = (a_{\ell})_{\sum_{p=0}^{\ell} h_{p}}(a_{\ell +1})_{h_{\ell+1}} \cdots (a_{m-1})_{h_{m-1}}
 (a_m)_{\sum_{p=m}^{k-1} h_{p}}.$$

 Now, it follows
 $$\partial_{{m-\ell}}^{k-1}(\alpha) + \partial_{{m-\ell}}^{k-1}(\beta)  = (a_{\ell})_{h'}(a_{\ell +1})_{h_{\ell+1}} \cdots (a_{m-1})_{h_{m-1}}
 (a_m)_{h''},$$
  where
 $$h' = \min\left\{\sum_{p=0}^{\ell} i_p \, , \,\sum_{p=0}^{\ell} j_p\right\} \; \mbox{and} \; h'' = \min\left\{\sum_{p=m}^{k-1} i_p \, , \,\sum_{p=m}^{k-1} j_p\right\}
 $$
 and $h_{\ell+1}, \ldots, h_{m-1}$ are well defined. Since
 $$\min\left\{\sum_{p=0}^{\ell} i_p \, , \,\sum_{p=0}^{\ell} j_p\right\} = \sum_{p=0}^{\ell} h_{p} \;\; \mbox{and} \;\;
 \min\left\{\sum_{p=m}^{k-1} i_p \, , \,\sum_{p=m}^{k-1} j_p\right\} = \sum_{p=m}^{k-1} h_{p},$$
 it follows $\partial_{{m-\ell}}^{k-1}(\alpha) + \partial_{{m-\ell}}^{k-1}(\beta) = \partial_{{m-\ell}}^{k-1}(\alpha + \beta)$.

\vspace{4mm}

Let us consider the set below:
$$S_{{m-\ell}}^{k-1} = \left\{ \alpha \left|\; \alpha \in \sigma^{(n)}\{a_0,\ldots,a_{k-1}\}, \alpha (a_0) = \cdots = \alpha(a_{\ell}) \geq a_{\ell+1}, \alpha(a_m) \leq a_m \right.\right\}.$$

For $\ds \alpha, \beta \in S_{{m-\ell}}^{k-1}$ we find
$$(\alpha + \beta)(a_{p}) = \alpha(a_{p}) + \beta(a_{p}) = \alpha(a_{q}) + \beta(a_{q}) = (\alpha + \beta)(a_{q}),$$
where $0 \leq p \leq \ell$, $0 \leq q \leq \ell$ and
$$(\alpha + \beta)(a_{p}) = \alpha(a_{p}) + \beta(a_{p}) \geq a_{\ell+1}\; \mbox{for all}\; p = 0, \ldots, \ell.$$
Moreover
$$(\alpha + \beta)(a_{m}) = \alpha(a_{m}) + \beta(a_{{m}}) \leq a_{m}.$$

Also it follows
$$(\alpha \beta)(a_{p}) = \beta(\alpha(a_{p})) = \beta(\alpha(a_{q})) = (\alpha \beta)(a_{q}),$$
where $0 \leq p \leq \ell$, $0 \leq q \leq \ell$ and
$$(\alpha \beta)(a_{p}) = \beta(\alpha(a_{p})) \geq  \beta(a_{\ell+1}) \geq \beta(a_p) \geq a_{\ell+1}\; \mbox{for all}\; p = 0, \ldots,  \ell.$$
Moreover
$$(\alpha \beta)(a_{m}) = \beta(\alpha(a_{m})) \leq \beta(a_m) \leq a_m.$$

So, we prove

\vspace{4mm}

\textbf{Lemma 3.} \textsl{The set $S_{{m-\ell}}^{k-1}$ is a subsemiring of the semiring $\sigma^{(n)}\{a_0,\ldots,a_{k-1}\}$.  }

\vspace{3mm}

Let us consider the set below:
$$R_{{m -\ell}}^{k-1} = \left\{ \alpha \left|\; \alpha \in \sigma^{(n)}\{a_0,\ldots,a_{k-1}\}, \alpha(a_{\ell}) \leq a_{\ell}, \; \alpha(a_{m}) \leq a_{m} \right.\right\}.$$

For $\ds \alpha, \beta \in R_{{m-\ell}}^{k-1}$ we find
$$(\alpha + \beta)(a_{\ell}) = \alpha(a_{\ell}) + \beta(a_{{\ell}}) \leq a_{\ell} \; \mbox{and analogously}\;
(\alpha + \beta)(a_{m})  \leq a_{m},$$
$$(\alpha \beta)(a_{\ell}) = \beta(\alpha(a_{\ell})) \leq  \beta(a_{\ell}) \leq a_{\ell}  \; \mbox{and analogously}\; (\alpha \beta)(a_{m})  \leq a_{m}.$$

Thus, we prove

\vspace{5mm}

\textbf{Lemma 4.} \textsl{The set $R_{{m-\ell}}^{k-1}$ is a subsemiring of the semiring $\sigma^{(n)}\{a_0,\ldots,a_{k-1}\}$.  }

\vspace{5mm}

Let $\mathcal{D}_{{m-\ell}}^{k-1} = S_{{m-\ell}}^{k-1}\cup R_{{m-\ell}}^{k-1}$. Note that $S_{{m-\ell}}^{k-1}\cap R_{{m-\ell}}^{k-1} = \varnothing$.

\vspace{5mm}

\textbf{Lemma 5.} \textsl{The set $\mathcal{D}_{{m-\ell}}^{k-1}$ is a subsemiring of the semiring $\sigma^{(n)}\{a_0,\ldots,a_{k-1}\}$.  }

\vspace{2mm}

\emph{Proof.} Let $\alpha \in S_{{m-\ell}}^{k-1}$ and $\beta \in R_{{m-\ell}}^{k-1}$.

Since $(\alpha + \beta)(a_{p}) = \alpha(a_{p}) + \beta(a_{p})$, $\alpha(a_p) \geq a_{\ell+1}$ and $\beta(a_p) \leq \beta(a_{\ell}) \leq a_{\ell}$, it follows $(\alpha + \beta)(a_{p}) = \alpha(a_{p}) \geq a_{\ell+1}$. As $\alpha(a_p) = \alpha(a_q)$ we have $(\alpha + \beta)(a_{p}) = (\alpha + \beta)(a_{q})$, where $0 \leq p \leq \ell$, $0 \leq q \leq \ell$. Moreover,
$(\alpha + \beta)(a_{m}) = \alpha(a_{m}) + \beta(a_{{m}}) \leq a_{m}$.
 Hence $\alpha + \beta \in S_{{m -\ell}}^{k-1}$.

Let $0 \leq p \leq \ell$. As $\beta(a_p) \leq \beta(a_{\ell}) \leq a_{\ell}$, it follows $\beta(a_p) = a_s$, where $s \leq \ell$. Analogously, for any $q$, $0 \leq q \leq \ell$, it follows $\beta(a_q) = a_t$, where $t \leq  \ell$. Then
$$(\beta\alpha)(a_p) =  \alpha(\beta(a_p)) = \alpha(a_s) = \alpha(a_t) =  \alpha(\beta(a_q)) = (\beta\alpha)(a_q).$$
On the other hand, $(\beta\alpha)(a_p) =  \alpha(\beta(a_p)) = \alpha(a_s) \geq a_{\ell+1}$.  Moreover,
$(\beta\alpha )(a_{m}) = \alpha(\beta(a_{m})) \leq \alpha(a_{{m}}) \leq a_{m}$. Hence $\beta\alpha \in S_{{m - \ell}}^{k-1}$.

Let $\alpha(a_p) = a_{\ell+1}$, where $0 \leq p \leq \ell$, and $\beta(a_{\ell+1}) \geq a_{\ell+1}$. Then
$$(\alpha\beta)(a_p) =  \beta(\alpha(a_p)) = \beta(a_{\ell+1}) \geq a_{\ell+1}.$$

 As $(\alpha\beta)(a_q) =  \beta(\alpha(a_q)) = \beta(a_{\ell+1})$\ where $0 \leq q \leq \ell$, it follows $(\alpha\beta)(a_p) = (\alpha\beta)(a_q)$. Moreover,
$(\alpha\beta)(a_{m}) = \beta(\alpha(a_{m})) \leq \beta(a_{{m}}) \leq a_{m}$.
Hence $\alpha\beta \in S_{{m - \ell}}^{k-1}$.

Let $\alpha(a_p) = a_{\ell+1}$, where $0 \leq p \leq \ell$, and $\beta(a_{\ell+1}) \leq a_{\ell}$. Then
$$(\alpha\beta)(a_{\ell}) =  \beta(\alpha(a_{\ell})) = \beta(a_{\ell+1}) \leq a_{\ell} \; \mbox{and}$$
$$(\alpha\beta)(a_{m}) =  \beta(\alpha(a_{m})) \leq \beta(a_{m}) \leq a_{m}.$$

Hence $\alpha\beta \in R_{m-\ell}^{k-1}$.

Let $\alpha(a_p) = a_{r}$, where $r > \ell + 1$ and $0 \leq p \leq \ell$. Let  $\beta(a_{r}) \geq a_{\ell+1}$. Then
$$(\alpha\beta)(a_p) =  \beta(\alpha(a_p)) = \beta(a_{r}) \geq a_{\ell+1}.$$

 As $(\alpha\beta)(a_q) =  \beta(\alpha(a_q)) = \beta(a_{r})$, it follows $(\alpha\beta)(a_p) = (\alpha\beta)(a_q)$, where\\ $0 \leq p \leq \ell$, $0 \leq q \leq \ell$. Moreover,
$(\alpha\beta)(a_{m}) = \beta(\alpha(a_{m})) \leq \beta(a_{{m}}) \leq a_{m}$.
Hence $\alpha\beta \in S_{{m - \ell}}^{k-1}$.

Let $\alpha(a_p) = a_{r}$, where $r > \ell + 1$ and $0 \leq p \leq \ell$. Let  $\beta(a_{r}) \leq a_{\ell}$. Then
$$(\alpha\beta)(a_{\ell}) =  \beta(\alpha(a_{\ell})) = \beta(a_{r}) \leq a_{\ell}.$$
 Moreover,
$(\alpha\beta)(a_{m}) = \beta(\alpha(a_{m})) \leq \beta(a_{{m}}) \leq a_{m}$.
Hence $\alpha\beta \in R_{{m-\ell}}^{k-1}$.

\vspace{4mm}

\textbf{Lemma 6.} \textsl{For any endomorphisms $\alpha, \beta \in \mathcal{D}_{{{m-\ell}}}^{k-1}$, it follows}
 $$\partial_{{m-\ell}}^{k-1}(\alpha \beta) = \partial_{{m-\ell}}^{k-1}(\alpha)\beta + \alpha\partial_{{m-\ell}}^{k-1}(\beta), $$
where $0 \leq \ell < m \leq k-1$.

\vspace{2mm}

\emph{Proof.} Let $\alpha = (a_0)_{i_0} \cdots (a_{k-1})_{i_{k-1}}$ and $\beta = (a_0)_{j_0} \cdots (a_{k-1})_{j_{k-1}}$, where\\ $\ds \sum_{p=0}^{k-1} i_p = \sum_{p=0}^{k-1} j_p = n$.

\vspace{2mm}

\emph{Case 1.} Let $\beta \in S_{{m-\ell}}^{k-1}$ and let assume that $\beta(a_0) = \cdots = \beta(a_{\ell}) = \cdots = \beta(a_{q_0}) = a_q$, where $q_0 \geq \ell$ be the greatest number with this property and $q \geq \ell+1$.

\vspace{2mm}

\emph{Case 1.1.} Let us assume that  the image of any element of the set $\{a_{m+1}, \ldots, a_{k-1}\}$ under the endomorphism $\beta$ is smaller than or equal to $a_{m}$, so, we suppose that  $\beta(a_{k-1}) = a'_{k-1} = a_r \leq a_{m}$.

 Since $\ds  \sum_{p=0}^{q-1} j_{p} \leq a_0, \ldots, a_{q_0} \leq \sum_{p=0}^{q} j_{p} - 1$,
it follows
$$\alpha\beta = (a_q)_{\sum_{p=0}^{q_0} i_{p}}(a_{q + 1})_{s_{q+1}} \cdots (a_{r})_{s_{r}},$$
where
$\ds \sum_{p=0}^{q_0} i_{p} + \sum_{p=q+1}^r s_p = n$. Now, it follows $\partial_{{m-\ell}}^{k-1}(\alpha \beta) = \alpha \beta$. Since
 $$\partial_{{m-\ell}}^{k-1}(\alpha) = (a_{\ell})_{\sum_{p=0}^{\ell} i_{p}}(a_{\ell +1})_{i_{\ell+1}} \cdots (a_{m-1})_{i_{m-1}}
 (a_m)_{\sum_{p=m}^{k-1} i_{p}},$$ using that $\ell \leq q_0$ and $r \leq m$,
 it follows $\partial_{{m-\ell}}^{k-1}(\alpha)\beta = \alpha\beta$. It is clear that
 $$\partial_{{m-\ell}}^{k-1}(\beta) = (a_{\ell})_{\sum_{p=0}^{\ell} j_{p}}(a_{\ell +1})_{j_{\ell+1}} \cdots (a_{m-1})_{
  j_{m-1}}  (a_m)_{\sum_{p=m}^{k-1} j_{ p}}$$
 and since $\ds a_0 \geq \sum_{p=0}^{q-1} j_{p} \geq \sum_{p=0}^{\ell} j_{p}$, it follows $\alpha\partial_{{m-\ell}}^{k-1}(\beta) = \alpha\beta$.

 Hence
$\partial_{{m-\ell}}^{k-1}(\alpha \beta) = \partial_{{m-\ell}}^{k-1}(\alpha)\beta + \alpha\partial_{{m-\ell}}^{k-1}(\beta)$.

\vspace{2mm}

\emph{Case 1.2.}
Let us assume that the image of any element of the set $\{a_{m+1}, \ldots, a_{k-1}\}$ under the endomorphism $\beta$ is greater   than $a_{m}$. Then
$$\alpha\beta = (a_q)_{\sum_{p=0}^{q_0} i_{p}}(a_{q + 1})_{s_{q+1}} \cdots (a_{k-1})_{s_{k-1}},$$

where
$\ds \sum_{p=0}^{q_0} i_{p} + \sum_{p=q+1}^{k-1} s_p = n$ Hence
$$\partial_{{m-\ell}}^{k-1}(\alpha \beta) = (a_q)_{\sum_{p=0}^{q_0} i_{p}}(a_{q + 1})_{s_{q+1}} \cdots (a_{m})_{\sum_{p=m}^{k-1} s_{p}}.$$

Using that $\ds \partial_{{m-\ell}}^{k-1}(\alpha) = (a_{\ell})_{\sum_{p=0}^{\ell} i_{p}}(a_{\ell +1})_{i_{\ell+1}} \cdots (a_{m-1})_{i_{m-1}}
 (a_m)_{\sum_{p=m}^{k-1} i_{p}}$ and $\ell \leq q_0$, it follows
$$\partial_{{m-\ell}}^{k-1}(\alpha) \beta = (a_q)_{\sum_{p=0}^{q_0} i_{p}}(a_{q + 1})_{s_{q+1}} \cdots (a_{m})_{\sum_{p=m}^{k-1} s_{p}}.$$

  Note  that the indices $\ds s_{q+1}, \ldots, \sum_{p=m}^{k-1} s_{p}$ are the same as in $\partial_{{m-\ell}}^{k-1}(\alpha\beta)$ because $\partial_{{m-\ell}}^{k-1}(\alpha)(a_p) = \alpha(a_p)$ for $p = q+1, \ldots, m-1$.

Since
 $$\partial_{{m-\ell}}^{k-1}(\beta) = (a_{\ell})_{\sum_{p=0}^{\ell} j_{p}}(a_{\ell +1})_{j_{\ell+1}} \cdots (a_{m-1})_{
  j_{m-1}}  (a_m)_{\sum_{p=m}^{k-1} j_{ p}}$$ and using that  $\partial_{{m-\ell}}^{k-1}(\beta)(a_p) = \beta(a_p)$ for $p = q+1, \ldots, m-1$, it follows
$$\alpha\partial_{{m-\ell}}^{k-1}(\beta) = (a_q)_{\sum_{p=0}^{q_0} i_{p}}(a_{q + 1})_{s_{q+1}} \cdots (a_{m})_{\sum_{p=m}^{k-1} s_{p}}.$$

 Hence
$\partial_{{m-\ell}}^{k-1}(\alpha \beta) = \partial_{{m-\ell}}^{k-1}(\alpha)\beta + \alpha\partial_{{m-\ell}}^{k-1}(\beta)$.

\vspace{3mm}

\emph{Case 1.3.}  Let us assume that only images of  some of  the elements of  $\{a_{m+1}, \ldots, a_{k-1}\}$ under the endomorphism $\beta$ are smaller than or equal to $a_{m}$. So, we suppose that  $\beta(a_{k-1}) = a'_{k-1} = a_r > a_{m}$ and there exist a number $t$, $m + 1 \leq t < k-1$, such that $\beta(a_{t}) < a_{m}$. Let $t$ be the smallest number with this property.  Hence,  $a'_{t} < a_{m} \leq a'_{t+1}$.

Now, considering (3), it follows
$$\alpha\beta = (a'_{0})_{\sum_{p=0}^{n_0} i_p} \cdots  (a'_{u})_{\sum_{p=n_{u-1}+1}^{n_u} i_p} \cdots  (a'_{k-1})_{\sum_{p=n_{r-1}+1}^{k-1} i_p},$$
where $a'_0 = a_q$, $n_0 = q_0$ and $1 \leq u \leq k-1$. Then
$$\partial_{{m-\ell}}^{k-1}(\alpha \beta) = (a_q)_{\sum_{p=0}^{n_0} i_{p}}(a'_{q_0 + 1})_{\sum_{p=n_0+1}^{n_1} i_p} \cdots (a'_{t})_{\sum_{p=n_{t-1}+1}^{n_{t}} i_p} (a_{m})_{\sum_{p=n_{t} + 1}^{k-1} i_{p}}.$$

Since $\ds \partial_{{m-\ell}}^{k-1}(\alpha) = (a_{\ell})_{\sum_{p=0}^{\ell} i_{p}}(a_{\ell +1})_{i_{\ell+1}} \cdots
 (a_m)_{\sum_{p=m}^{k-1} i_{p}}$, $\;\ell \leq q_0 = n_0$ and $\ds \partial_{{m-\ell}}^{k-1}(\alpha)(a_u) = \alpha(a_u)$ for  $u = n_0 + 1, \ldots, n_{t}$, it follows
$$\partial_{{m-\ell}}^{k-1}(\alpha)\, \beta = (a_q)_{\sum_{p=0}^{n_0} i_{p}}(a'_{q_0+1})_{\sum_{p=n_0+1}^{n_1} i_p} \cdots (a'_{t})_{\sum_{p=n_{t-1}+1}^{n_{t}} i_p} (a_{m})_{\sum_{p=n_{t} + 1}^{k-1} i_{p}}.$$

Since
$\ds \partial_{{m-\ell}}^{k-1}(\beta) = (a_{\ell})_{\sum_{p=0}^{\ell} j_{p}}(a_{\ell +1})_{j_{\ell+1}} \cdots (a_{m-1})_{
  j_{m-1}}  (a_m)_{\sum_{p=m}^{k-1} j_{ p}}$,\\  $\;\ell \leq q_0 = n_0$ and $\ds \partial_{{m-\ell}}^{k-1}(\beta)(a_u) = \beta(a_u)$ for  $u = \ell + 1, \ldots, m-1$, it follows
$$\alpha \,\partial_{{m-\ell}}^{k-1}(\beta)= (a_q)_{\sum_{p=0}^{n_0} i_{p}}(a'_{q_0+1})_{\sum_{p=n_0+1}^{n_1} i_p} \cdots (a'_{t})_{\sum_{p=n_{t-1}+1}^{n_{t}} i_p} (a_{m})_{\sum_{p=n_{t} + 1}^{k-1} i_{p}}.$$

 Hence
$\partial_{{m-\ell}}^{k-1}(\alpha \beta) = \partial_{{m-\ell}}^{k-1}(\alpha)\beta + \alpha\partial_{{m-\ell}}^{k-1}(\beta)$.

\vspace{4mm}

\emph{Case 2.} Let $\beta \in R_{{m-\ell}}^{k-1}$.

\vspace{2mm}

\emph{Case 2.1.} Let us assume that the image of any element of the set $\{a_{m+1}, \ldots, a_{k-1}\}$ under the endomorphism $\beta$ is smaller than  or equal to $a_{\ell}$. So, we suppose  that $\beta(a_{k-1}) = a'_{k-1} = a_r$, where $r \leq \ell$.  Thus $\alpha\beta = (a_0)_{s_0} \cdots (a_r)_{s_{r}}$ which implies $\partial_{{m-\ell}}^{k-1}(\alpha \beta) = (a_{\ell})_n$.
It  is clear that
$$ \partial_{{m-\ell}}^{k-1}(\alpha) = (a_{\ell})_{\sum_{p=0}^{\ell} i_{p}}(a_{\ell +1})_{i_{\ell+1}} \cdots
 (a_m)_{\sum_{p=m}^{k-1} i_{p}}$$
 and then $\partial_{{m-\ell}}^{k-1}(\alpha) \beta \leq (a_{\ell})_n$. Since
$$ \partial_{{m-\ell}}^{k-1}(\beta) = (a_{\ell})_{\sum_{p=0}^{\ell} j_{p}}(a_{\ell +1})_{j_{\ell+1}} \cdots (a_{m-1})_{
  j_{m-1}}  (a_m)_{\sum_{p=m}^{k-1} j_{ p}}$$
  and  $\ds a_{k-1} \leq \sum_{p=0}^{\ell} j_{p} - 1$, it follows
  $\alpha\,\partial_{{m-\ell}}^{k-1}(\beta) = (a_{\ell})_n$.
\vspace{1mm}

 Hence
$\partial_{{m-\ell}}^{k-1}(\alpha \beta) = \partial_{{m-\ell}}^{k-1}(\alpha)\beta + \alpha\partial_{{m-\ell}}^{k-1}(\beta)$.

\vspace{2mm}

\emph{Case 2.2.} Let us assume that the image of any element of the set $\{a_{m+1}, \ldots, a_{k-1}\}$ under the endomorphism $\beta$ is smaller than or equal to $a_m$. So, we suppose   that $\beta(a_{k-1}) = a'_{k-1} = a_r$, where $ \ell < r \leq m$.

Let  $\beta(a_t) \leq a_\ell$, where $\ell \leq t < k-1$ and $t$ be the greatest number with this property. Hence, $a'_t \leq a_\ell < a'_{t+1}$. Then, from (3), it follows
$$\alpha\beta = (a'_{0})_{\sum_{p=0}^{n_0} i_p} \cdots  (a'_{t})_{\sum_{p=n_{t-1}+1}^{n_t} i_p} \cdots
 (a'_{k-1})_{\sum_{p=n_{r-1}+1}^{k-1} i_p},$$
 which implies
$$ \partial_{{m-\ell}}^{k-1}(\alpha\beta) = (a_{\ell})_{\sum_{p=0}^{n_t} i_{p}} \cdots  (a_r)_{\sum_{p=n_{r-1}+1}^{k-1} i_p}.$$

Clearly $\ds \partial_{{m-\ell}}^{k-1}(\alpha) = (a_{\ell})_{\sum_{p=0}^{\ell} i_{p}}(a_{\ell +1})_{i_{\ell+1}} \cdots
 (a_m)_{\sum_{p=m}^{k-1} i_{p}}$, but, if $\beta(a_\ell)  = a_{\ell_0} < a_\ell$, it follows
$$ \partial_{{m-\ell}}^{k-1}(\alpha)\beta = (a_{\ell_0})_{\sum_{p=0}^{\ell} i_{p}} \cdots  (a'_m)_{\sum_{p=m}^{k-1} i_p} < \partial_{{m-\ell}}^{k-1}(\alpha\beta).$$
Note, that generally we observe $\partial_{{m-\ell}}^{k-1}(\alpha)\beta \leq \partial_{{m-\ell}}^{k-1}(\alpha\beta)$.

Since
$$ \partial_{{m-\ell}}^{k-1}(\beta) = (a_{\ell})_{\sum_{p=0}^{\ell} j_{p}}(a_{\ell +1})_{j_{\ell+1}} \cdots (a_{m-1})_{
  j_{m-1}}  (a_m)_{\sum_{p=m}^{k-1} j_{ p}}$$
and using that $a'_t \leq a_\ell$ and $a'_{k-1} = a_r \leq a_m$, it follows
$$\alpha\, \partial_{{m-\ell}}^{k-1}(\beta) = (a_{\ell})_{\sum_{p=0}^{n_t} i_{p}} \cdots  (a_r)_{\sum_{p=n_{r-1}+1}^{k-1} i_p}.$$

 Hence
$\partial_{{m-\ell}}^{k-1}(\alpha \beta) = \partial_{{m-\ell}}^{k-1}(\alpha)\beta + \alpha\partial_{{m-\ell}}^{k-1}(\beta)$.

\vspace{2mm}

\emph{Case 2.3.}
Let us assume that the image of any element of the set $\{a_{m+1}, \ldots, a_{k-1}\}$ under the endomorphism $\beta$ is greater than or equal to $a_{m}$. Then, considering reasonings similar to those in Case 2.2 and (3), it follows
$$\alpha\beta = (a'_{0})_{\sum_{p=0}^{n_0} i_p} \cdots  (a'_{t})_{\sum_{p=n_{t-1}+1}^{n_t} i_p} \cdots  (a'_{k-1})_{\sum_{p=n_{r-1}+1}^{k-1} i_p},$$
where $\ell \leq t \leq m$ and $t$ be the smallest number such that  $\beta(a_t) \leq a_\ell$.

 Now, it follows
$$\partial_{{m-\ell}}^{k-1}(\alpha \beta) = (a_\ell)_{\sum_{p=0}^{n_t} i_{p}}(a'_{t + 1})_{\sum_{p=n_t+1}^{n_{t+1}} i_p} \cdots (a'_{m-1})_{\sum_{p=n_{m-2}+1}^{n_{m-1}} i_p} (a_{m})_{\sum_{p=n_{m-1} + 1}^{k-1} i_{p}}.$$

Since $ \partial_{{m-\ell}}^{k-1}(\alpha) = (a_{\ell})_{\sum_{p=0}^{\ell} i_{p}}(a_{\ell +1})_{i_{\ell+1}} \cdots
 (a_m)_{\sum_{p=m}^{k-1} i_{p}}$, as in the previous case we find $\partial_{{m-\ell}}^{k-1}(\alpha)\beta \leq \partial_{{m-\ell}}^{k-1}(\alpha\beta)$.

Since
$ \partial_{{m-\ell}}^{k-1}(\beta) = (a_{\ell})_{\sum_{p=0}^{\ell} j_{p}}(a_{\ell +1})_{j_{\ell+1}} \cdots (a_{m-1})_{
  j_{m-1}}  (a_m)_{\sum_{p=m}^{k-1} j_{ p}}$
and using that $a'_t \leq a_\ell$ and $a'_{m+1} \geq a_m$, it follows
$$\alpha\, \partial_{{m-\ell}}^{k-1}(\beta) = (a_\ell)_{\sum_{p=0}^{n_t} i_{p}}(a'_{t + 1})_{\sum_{p=n_t+1}^{n_{t+1}} i_p}\!\! \cdots (a'_{m-1})_{\sum_{p=n_{m-2}+1}^{n_{m-1}} i_p} (a_{m})_{\sum_{p=n_{m-1} + 1}^{k-1} i_{p}}.$$

 Hence
$\partial_{{m-\ell}}^{k-1}(\alpha \beta) = \partial_{{m-\ell}}^{k-1}(\alpha)\beta + \alpha\partial_{{m-\ell}}^{k-1}(\beta)$.

\vspace{2mm}

\emph{Case 2.4.} Let us assume that only  images of  some of  the elements of  $\{a_{m+1}, \ldots, a_{k-1}\}$ under the endomorphism $\beta$ are smaller than  or equal to $a_{m}$. So, we suppose that  $\beta(a_{k-1}) = a'_{k-1} = a_r > a_{m}$ and there exist a number $t$, $m + 1 \leq t < k-1$, such that $\beta(a_{t}) < a_{m}$. Let $t$ be the smallest number with this property.  Hence,  $a'_{t} < a_{m} \leq a'_{t+1}$. Let $q$ be the greatest number such that $\beta(a_q) \leq a_\ell$, where $\ell \leq q < t$. Hence, $a'_q \leq a_\ell < a'_{q+1}$.
Then, from (3), it follows
$$\alpha\beta = (a'_{0})_{\sum_{p=0}^{n_0} i_p} \cdots  (a'_{q})_{\sum_{p=n_{q-1}+1}^{n_q} i_p} \cdots (a'_{t})_{\sum_{p=n_{t-1}+1}^{n_t} i_p} \cdots  (a'_{k-1})_{\sum_{p=n_{r-1}+1}^{k-1} i_p}.$$

Now, it follows
$$\partial_{{m-\ell}}^{k-1}(\alpha \beta) = (a_\ell)_{\sum_{p=0}^{n_q} i_{p}}(a'_{q + 1})_{\sum_{p=n_q+1}^{n_{q+1}} i_p} \cdots (a'_{t})_{\sum_{p=n_{t-1}+1}^{n_{t}} i_p} (a_{m})_{\sum_{p=n_{t} + 1}^{k-1} i_{p}}.$$

Since $\ds \partial_{{m-\ell}}^{k-1}(\alpha) = (a_{\ell})_{\sum_{p=0}^{\ell} i_{p}}(a_{\ell +1})_{i_{\ell+1}} \cdots
 (a_m)_{\sum_{p=m}^{k-1} i_{p}}$, if $\beta(a_\ell) = a_{\ell_0} < a_\ell$, it follows
$$\partial_{{m-\ell}}^{k-1}(\alpha)\, \beta = (a_{\ell_0})_{\sum_{p=0}^{\ell} i_p} \cdots (a'_{t})_{\sum_{p=n_{t-1}+1}^{n_{t}} i_p} (a_{m})_{\sum_{p=n_{t} + 1}^{k-1} i_{p}} < \partial_{{m-\ell}}^{k-1}(\alpha\beta).$$
Generally we observe $\partial_{{m-\ell}}^{k-1}(\alpha)\beta \leq \partial_{{m-\ell}}^{k-1}(\alpha\beta)$.

It is clear that
$$ \partial_{{m-\ell}}^{k-1}(\beta) = (a_{\ell})_{\sum_{p=0}^{\ell} j_{p}}(a_{\ell +1})_{j_{\ell+1}} \cdots (a_{m-1})_{
  j_{m-1}}  (a_m)_{\sum_{p=m}^{k-1} j_{ p}}.$$

First multiplier in $\partial_{{m-\ell}}^{k-1}(\alpha \beta)$ in its representation as a product is $(a_\ell)_{\sum_{p=0}^{n_q} i_{p}}$. Since
$$\sum_{p=0}^{n_q} i_{p} = \sum_{p=0}^{n_0} i_{p} + \cdots + \sum_{p=n_{q-1} +1}^{n_q} i_{p},$$
it means that $\partial_{{m-\ell}}^{k-1}(\alpha \beta)(a'_p) = a_\ell$, where $p = 0, \ldots, q$. Since, for $p = 0, \ldots, q$, we have $\ds a'_p < \sum_{p=0}^\ell j_p$, it follows $\alpha\partial_{{m-\ell}}^{k-1}(\beta)(a'_p) = a_\ell$. So, first multiplier in the representation as a product of $\alpha\partial_{{m-\ell}}^{k-1}(\beta)$ is $(a_\ell)_{\sum_{p=0}^{n_q} i_{p}}$.

On the other hand $\ds \partial_{{m-\ell}}^{k-1}(\beta)(a_u) = \beta(a_u)$ for  $u = q + 1, \ldots, m = 1$. Thus, it follows
$$\alpha \,\partial_{{m-\ell}}^{k-1}(\beta)= (a_q)_{\sum_{p=0}^{n_0} i_{p}}(a'_{q_0+1})_{\sum_{p=n_0+1}^{n_1} i_p} \cdots (a'_{t})_{\sum_{p=n_{t-1}+1}^{n_{t}} i_p} (a_{m})_{\sum_{p=n_{t} + 1}^{k-1} i_{p}}.$$

 Hence
$\partial_{{m-\ell}}^{k-1}(\alpha \beta) = \partial_{{m-\ell}}^{k-1}(\alpha)\beta + \alpha\partial_{{m-\ell}}^{k-1}(\beta)$.

\vspace{5mm}

\textbf{Theorem.}  \textsl{The map $\partial_{{m-\ell}}^{k-1} : \mathcal{D}_{{{m-\ell}}}^{k-1}  \rightarrow \sigma^{(n)}\{a_{\ell},\ldots,a_{m}\}$  is a derivation. The maximal subsemiring of
 $\sigma^{(n)}\{a_0, a_1, \ldots, a_{k-1}\}$ closed under the derivation $\partial_{{m-\ell}}^{k-1}$ is a semiring $\mathcal{D}_{{{m-\ell}}}^{k-1}$.}

\vspace{2mm}

\emph{ {{Proof.}}} From the lemmas  we obtain that  $\partial_{{m-\ell}}^{k-1} : \mathcal{D}_{{{m-\ell}}}^{k-1}\!\!  \rightarrow\! \sigma^{(n)}\{a_{\ell},\ldots,a_{m}\}$ is a derivation.

 If we suppose that $\beta \notin \mathcal{D}_{{{m-\ell}}}^{k-1} = S_{{m-\ell}}^{k-1}\cup R_{{m-\ell}}^{k-1}$, then it follows $\beta(a_{m}) > a_{m}$. Let $t > m$ be the smallest number such that $\beta(a_m) = a_t$, so $a'_{t-1} \leq a_m < a_t$.   Let  $q$ be the greatest number such that $\beta(a_q) \leq a_\ell$, where $\ell \leq q$. Hence, $a'_q \leq a_\ell < a'_{q+1}$. From (3) we have
 $$\alpha\beta = (a'_{0})_{\sum_{p=0}^{n_0} i_p} \cdots  (a'_{q})_{\sum_{p=n_{q-1}+1}^{n_q} i_p} \cdots (a'_{t})_{\sum_{p=n_{t-1}+1}^{n_t} i_p} \cdots  (a'_{k-1})_{\sum_{p=n_{r-1}+1}^{k-1} i_p}.$$

Hence, it follows
$$\partial_{{m-\ell}}^{k-1}(\alpha \beta) = (a_\ell)_{\sum_{p=0}^{n_q} i_{p}}(a'_{q + 1})_{\sum_{p=n_q+1}^{n_{q+1}} i_p} \cdots  (a_{m})_{\sum_{p=n_{t-1} + 1}^{k-1} i_{p}}.$$

Since
$$ \partial_{{m-\ell}}^{k-1}(\alpha) = (a_{\ell})_{\sum_{p=0}^{\ell} i_{p}}(a_{\ell +1})_{i_{\ell+1}} \cdots
 (a_m)_{\sum_{p=m}^{k-1} i_{p}},$$
it follows by reasonings similar to those from Case 2.4 of Lemma 6, that
$$\partial_{{m-\ell}}^{k-1}(\alpha)\, \beta = (a_\ell)_{\sum_{p=0}^{n_q} i_{p}} \cdots  (a_{t})_{\sum_{p=n_{t-1} + 1}^{k-1} i_{p}},$$
when $\beta(a_\ell) = a_\ell$, or
$$\partial_{{m-\ell}}^{k-1}(\alpha)\, \beta = (a_{\ell_0})_{\sum_{p=0}^{\ell} i_{p}} \cdots  (a_{t})_{\sum_{p=n_{t-1} + 1}^{k-1} i_{p}},$$
when $\beta(a_\ell) = a_{\ell_0} < a_\ell$. In both cases, since $a_t > a_m$ we have that $\partial_{{m-\ell}}^{k-1}(\alpha)\, \beta$ is not smaller than or equal to $\partial_{{m-\ell}}^{k-1}(\alpha \beta)$ and this completes the proof.

\vspace{2mm}

Immediately from the theorem follows that $\sigma^{(n)}\{a_{\ell},\ldots,a_{m}\} \subset \mathcal{D}_{{{m-\ell}}}^{k-1}$. But $\sigma^{(n)}\{a_{\ell},\ldots,a_{m}\}$ is a left ideal of $\sigma^{(n)}\{a_0,\ldots,a_{k-1}\}$, so it is a left ideal of semiring  $\mathcal{D}_{{{m-\ell}}}^{k-1}$. Analogously any subsimplex $\sigma^{(n)}\{a_{r},\ldots,a_{m}\}$, where $r = \ell + 1, \ldots, m-1$ is a left ideal of $\mathcal{D}_{{{m-\ell}}}^{k-1}$.

\vspace{5mm}

{\bf 4 \hspace{1mm} Projections onto  the smallest $m$ -- simpleces}

\vspace{3mm}

Here we consider the derivation from the previous section in the case when $\ell = 0$. So, we explore the map
$$\partial_{{m - 0}}^{k-1- 0} = \partial_{{m}}^{k-1} : \sigma^{(n)}\{a_0,\ldots,a_{k-1}\} \rightarrow \sigma^{(n)}\{a_{0},\ldots,a_{m}\},$$
where $0  < m \leq k-1$,
such that for any endomorphism\\  $\alpha \in  \sigma^{(n)}\{a_0,\ldots,a_{k-1}\}$,   $\alpha = (a_0)_{i_0} \cdots (a_{k-1})_{i_{k-1}}$, where $\ds \sum_{p=0}^{k-1} i_p = n$,
 $$\partial_{{m}}^{k-1}(\alpha) = (a_{0})_{i_{0}} \cdots (a_{m-1})_{i_{m-1}} (a_m)_{\sum_{p=m}^{k-1} i_{p}}.$$

Now, it follows  $$\mathcal{D}_{{m}}^{k-1} = \left\{ \alpha \left|\; \alpha \in \sigma^{(n)}\{a_0,\ldots,a_{k-1}\},  \alpha(a_{m}) \leq a_{m} \right.\right\}.$$

Immediately from the theorem, it follows that
$$\partial_{{m}}^{k-1} :  \mathcal{D}_{{m}}^{k-1}  \rightarrow \sigma^{(n)}\{a_{0},\ldots,a_{m}\}$$
is a derivation and $\mathcal{D}_{{m}}^{k-1}$  is the maximal subsemiring of
 $\sigma^{(n)}\{a_0, a_1, \ldots, a_{k-1}\}$ closed under this derivation.

 Analogously, for any integer $m_1$, where $0 < m_1 < m$, we consider the derivation
$$\partial_{{m_1}}^{k-1} :  \mathcal{D}_{{m_1}}^{k-1}  \rightarrow \sigma^{(n)}\{a_{0},\ldots,a_{m_1}\}.$$

 From theorem, using the last derivation we can construct
$$\partial_{{m_1-0}}^{m-0} = \partial_{{m_1}}^{m} :  \mathcal{D}_{{m_1}}^{m}  \rightarrow \sigma^{(n)}\{a_{0},\ldots,a_{m_1}\},$$
which is a derivation and $\mathcal{D}_{{m_1}}^{m} = \left\{ \alpha \left|\; \alpha \in \sigma^{(n)}\{a_0,\ldots,a_{m}\},  \alpha(a_{m_1}) \leq a_{m_1} \right.\right\}$ is the maximal subsemiring of
 $\sigma^{(n)}\{a_0, a_1, \ldots, a_{m}\}$ closed under this derivation.  So, we consider a new derivation
$$\partial_{{m_1}}^{k-1} :  \mathcal{D}_{{m}}^{k-1}\cap\mathcal{D}_{{m_1}}^{k-1}  \rightarrow \sigma^{(n)}\{a_{0},\ldots,a_{m_1}\}$$
 and obtain that for any $\alpha \in \mathcal{D}_{{m}}^{k-1}\cap \mathcal{D}_{{m_1}}^{k-1}$, it follows $\partial_{{m_1}}^{k-1}(\alpha) = \partial_{m_1}^{m}\left(\partial_{m}^{k-1}(\alpha)\right)$. Hence, $\partial_{{m_1}}^{k-1} = \partial_{m_1}^{m} \cdot\, \partial_{m}^{k-1}$.
Note that the last composition is possible  for arbitrary number $m$, where $m_1 < m < k-1$.

In the semiring $$ \mathcal{D}\! =\! \bigcap_{m=0}^{k-1}\! \mathcal{D}_{{m}}^{k-1}\! =\! \left\{ \alpha \left| \alpha \in \sigma^{(n)}\{a_0,\ldots,a_{k-1}\}, \alpha(a_{m}) \leq a_{m}, m = 0, \ldots, k - 1 \! \right.\right\}$$
we can represent any derivation as composition of two  or more derivations.

In the case when $k = n$, i.e.
$\sigma^{(n)}\{a_0,\ldots,a_{k-1}\} =\widehat{\mathcal{E}}_{\mathcal{C}_n}$ the semiring $\ds \bigcap_{m=0}^{n-1} \mathcal{D}_{{m}}^{n-1}$
is well known -- it is the subsemiring $\mathcal{ON}_n$ of over nilpotent endomorphisms, see [7]. Hence, semiring $\mathcal{ON}_n$ is closed under any derivation $\partial_{{m}}^{n-1}$, where $m = 1, \ldots , n - 2$.
The order of this semiring is the n-th Catalan number $\ds \frac{1}{n+1}\binom{2n}{n}$.

The semiring $\mathcal{ON}_n$ has a zero element $\overline{0}$ and we can consider the subsemiring $\mathcal{N}_n$ of nilpotent elements, i.e.
the endomorphisms $\alpha \in \widehat{\mathcal{E}}_{\mathcal{C}_n}$ such that $\alpha^m = \overline{0}$ for some positive integer $m$. It is known, [7], that $\mathcal{N}_n$ is an ideal of $\mathcal{ON}_n$ and it is consists of all endomorphisms $\alpha$ such that $\alpha(t) < t$ for any $t = 1, \ldots, n-1$. The order of the ideal $\mathcal{N}_n$ is the $(n-1)$-th Catalan number $\ds \frac{1}{n}\binom{2n - 2}{n-1}$.

\vspace{4mm}

\textbf{Corollary.}  \textsl{The ideal $\mathcal{N}_n$ is closed under arbitrary derivation $\partial_{m}^{n-1}$, where}
$m = 1, \ldots , n - 2$.
\vspace{2mm}

\emph{Proof.}  Let $\alpha = (a_0)_{i_0} \cdots (a_{n})_{i_{n}} \in \mathcal{N}_n$. Then
 $$\partial_{{m}}^{n-1}(\alpha) = (a_{0})_{i_{0}} \cdots (a_{m-1})_{i_{m-1}} (a_m)_{\sum_{p=m}^{n-1} i_{p}}.$$

For $\ds 0 \leq t \leq \sum_{p = 0}^{m-1} i_p -1$, it follows $\partial_{m}^{n-1}(\alpha) = \alpha$, so, $\left(\partial_{m}^{n-1}(\alpha)\right)(t) < t$.

For $\ds \sum_{p = 0}^{m-1} i_p   \leq t \leq n - 1$, it follows $\left(\partial_{m}^{n-1}(\alpha)\right)(t) = a_{n-m} \leq \alpha(t) < t$.

Hence $\ds \partial_{m}^{n-1}(\alpha) \in \mathcal{N}_n$.

\vspace{5mm}

{\bf 5 \hspace{1mm} Projections onto  the greatest $k - \ell$ -- simpleces}

\vspace{3mm}

Here we consider the derivation from the  section 3 in the case when $m = k-1$. We explore the map
$$\partial_{{k-1-\ell}}^{k-1- 0}  = \partial_{{k-\ell-1}}^{k-1} : \sigma^{(n)}\{a_0,\ldots,a_{k-1}\} \rightarrow \sigma^{(n)}\{a_{\ell},\ldots,a_{k-1}\},$$
where $0  \leq \ell < k-1$,
such that for any  $\alpha \in  \sigma^{(n)}\{a_0,\ldots,a_{k-1}\}$,\break   $\alpha = (a_0)_{i_0} \cdots (a_{k-1})_{i_{k-1}}$, where $\ds \sum_{p=0}^{k-1} i_p = n$,
 $$\partial_{{k-\ell -1}}^{k-1}(\alpha) = (a_{\ell})_{\sum_{p=0}^{\ell} i_{p}}(a_{\ell + 1})_{i_{\ell + 1}} \cdots (a_{k-1})_{i_{k-1}}.$$

So, we consider
$$S_{{k-\ell-1}}^{k-1} = \left\{ \alpha \left|\; \alpha \in \sigma^{(n)}\{a_0,\ldots,a_{k-1}\}, \alpha (a_0) = \cdots = \alpha(a_{\ell}) \geq a_{\ell+1} \right.\right\}, \eqno{(4)}$$
$$R_{{k -\ell - 1}}^{k-1} = \left\{ \alpha \left|\; \alpha \in \sigma^{(n)}\{a_0,\ldots,a_{k-1}\}, \alpha(a_{\ell}) \leq a_{\ell} \right.\right\} \eqno{(5)}$$ and $\mathcal{D}_{{k-\ell-1}}^{k-1} = S_{{k-\ell-1}}^{k-1}\cup R_{{k-\ell-1}}^{k-1}$. Immediately from (4) and (5) it follows
$$S_{{k-\ell_1-1}}^{k-1} \subset S_{{k-\ell -1}}^{k-1}, \eqno{(6)}$$ where $1 \leq  \ell < \ell_1 < k-1$ and
$$S_{{k-\ell_1 -1}}^{k-1} \cap R_{{k-\ell-1}}^{k-1} = \varnothing, \eqno{(7)}$$
for arbitrary $1 \leq \ell, \ell_1 < k-1$.

 Thus, using (6)  and (7),  we obtain
$$\mathcal{D}_{{{k-\ell-1}}}^{k-1}\cap \mathcal{D}_{{{k-\ell_1 -1}}}^{k-1} = S_{{k-\ell - 1}}^{k-1} \cup \left( R_{{k-\ell-1}}^{k-1}\cap R_{{k-\ell_1 - 1}}^{k-1}\right).$$

Immediately from the theorem it follows that
$$\partial_{{k-\ell-1}}^{k-1} : \mathcal{D}_{{{k-\ell-1}}}^{k-1}  \rightarrow \sigma^{(n)}\{a_{\ell},\ldots,a_{k-1}\}$$
is a derivation and $\mathcal{D}_{{k-\ell-1}}^{k-1}$  is the maximal subsemiring of
 $\sigma^{(n)}\{a_0, a_1, \ldots, a_{k-1}\}$ closed under this derivation.

 Analogously, for any integer $\ell_1$, where $\ell < \ell_1 < k-1$, we consider the derivation
$$\partial_{{k-\ell_1-1}}^{k-1} : \mathcal{D}_{{{k-\ell_1-1}}}^{k-1}  \rightarrow \sigma^{(n)}\{a_{\ell_1},\ldots,a_{k-1}\}.$$
Here $\mathcal{D}_{{k-\ell_1-1}}^{k-1} = S_{{k-\ell_1-1}}^{k-1}\cup R_{{k-\ell_1-1}}^{k-1}$, where
$$S_{{k-\ell_1-1}}^{k-1} = \left\{ \alpha \left|\; \alpha \in \sigma^{(n)}\{a_0,\ldots,a_{k-1}\}, \alpha (a_0) = \cdots = \alpha(a_{\ell_1}) \geq a_{\ell_1+1} \right.\right\},$$
$$R_{{k -\ell_1 - 1}}^{k-1} = \left\{ \alpha \left|\; \alpha \in \sigma^{(n)}\{a_0,\ldots,a_{k-1}\}, \alpha(a_{\ell_1}) \leq a_{\ell_1} \right.\right\}.$$

At last for arbitrary integers $\ell$ and $\ell_1$, where $\ell < \ell_1 < k-1$, using theorem, we can construct the derivation
$$\partial_{{k- 1 -\ell_1}}^{k - 1 - \ell} = \partial_{{k-\ell_1-1}}^{k - \ell - 1} : \mathcal{D}_{{{k-\ell_1-1}}}^{k -\ell -1}  \rightarrow \sigma^{(n)}\{a_{\ell_1},\ldots,a_{k-1}\}.$$
Here $\mathcal{D}_{{k-\ell_1-1}}^{k - \ell- 1} = S_{{k-\ell_1-1}}^{k - \ell -1}\cup R_{{k-\ell_1-1}}^{k - \ell - 1}$, where
$$S_{{k-\ell_1-1}}^{k -\ell-1} = \left\{ \alpha \left|\; \alpha \in \sigma^{(n)}\{a_\ell,\ldots,a_{k-1}\}, \alpha (a_\ell) = \cdots = \alpha(a_{\ell_1}) \geq a_{\ell_1+1} \right.\right\},$$
$$R_{{k -\ell_1 - 1}}^{k -\ell -1} = \left\{ \alpha \left|\; \alpha \in \sigma^{(n)}\{a_\ell,\ldots,a_{k-1}\}, \alpha(a_{\ell_1}) \leq a_{\ell_1} \right.\right\}.$$

For derivation
$$\partial_{{k-\ell_1-1}}^{k-1} : \mathcal{D}_{{{k-\ell-1}}}^{k-1}\cap \mathcal{D}_{{{k-\ell_1-1}}}^{k-1}  \rightarrow \sigma^{(n)}\{a_{\ell_1},\ldots,a_{k-1}\}$$
we obtain that for any endomorphism $\alpha \in \mathcal{D}_{{{k-\ell-1}}}^{k-1}\cap \mathcal{D}_{{{k-\ell_1-1}}}^{k-1}$, it follows
$\partial_{{k-\ell_1 -1}}^{k-1}(\alpha) = \partial_{{k-\ell_1-1}}^{k-\ell-1}\left(\partial_{{k-\ell-1}}^{k-1}(\alpha)\right)$.
Hence, $$\partial_{{k-\ell_1 -1}}^{k-1}  = \partial_{{k-\ell-1}}^{k-1} \, \cdot \, \partial_{{k-\ell_1-1}}^{k-\ell-1}$$
for any number $\ell_1$, where $\ell < \ell_1 \leq k-1$.

In the semiring $\ds \bigcap_{\ell=1}^{k-1} \mathcal{D}_{{{k-\ell-1}}}^{k-1}$
we can represent any derivation as composition  of two (or more) derivations.

From (6) we find
$$S_{{1}}^{k-1} \subset S_{{2}}^{k-1} \subset \cdots  \subset S_{{k-2}}^{k-1}.$$
Since
$$S_{{1}}^{k-1} = \left\{ \alpha \left| \alpha (a_0) = \cdots = \alpha(a_{k-2}) \geq a_{k-1} \right.\right\} = \{(a_{k-1})_n\}, \eqno{(8)}$$
it follows that $\ds \bigcap_{\ell=1}^{k-1} S_{{k-\ell-1}}^{k-1} =  \{(a_{k-1})_n\}$.

On the other hand, for $\mathcal{D}_{{m}}^{k-1}$ from section 4 and $R_{{k -\ell - 1}}^{k-1}$ from (5) we find that
$R_{{k -\ell - 1}}^{k-1} = \mathcal{D}_{{\ell}}^{k-1}$. Hence, using (7), it follows
$$ \bigcap_{\ell=1}^{k-1} \mathcal{D}_{{{k-\ell-1}}}^{k-1} = \bigcap_{\ell=1}^{k-1}\left(S_{{k-\ell-1}}^{k-1}\cup R_{{k-\ell-1}}^{k-1} \right) =$$ $$=  \{(a_{k-1})_n\} \cup \left(\bigcap_{\ell=1}^{k-1} R_{{k-\ell-1}}^{k-1} \right) = \{(a_{k-1})_n\} \cup \mathcal{D}.$$

So, as in the previous section, in the same semiring $\mathcal{D}$ we can represent any derivation of type $\partial_{{k-\ell-1}}^{k-1}$ as a composition of derivations.

In the case when $k = n$, the order of  $\mathcal{D} = \mathcal{ON}_n$ is the $n$-th Catalan number. The impotant role of the Catalan numbers appears also in the next reasonings.

For $k = n$ from (4) it follows
$$S_{{n-\ell-1}}^{n-1} = \left\{ \alpha \left|\; \alpha \in \sigma^{(n)}\{0,\ldots,n-1\}, \alpha (0) = \cdots = \alpha(\ell) \geq \ell+1 \right.\right\},$$
where $0 \leq \ell < n-1$. As in (8) we find $S_{{1}}^{n-1} = \{(n-1)_n\}$.

\vspace{5mm}

\textbf{Proposition.}  \textsl{For any $1 \leq p \leq n-2$ the order of semiring $S_{p}^{n-1}$ is}
$$\left|S_{{p}}^{n-1}\right| = p\,C_p,$$
\textsl{where $C_p$ is the $p$-th Catalan number.}

\vspace{5mm}

\begin{center}
  \textbf{References}
\end{center}

\vspace{2mm}

\z J. Golan, Semirings and their application, Kluwer Ac.Publ., Dodrecht, 1999.

\z J. J$\hat{\mbox{e}}$zek, T. Kepka, M.  Mar$\grave{\mbox{o}}$ti, The endomorphism semiring of a  semi\-lattice, Semigroup Forum, 78, 2009, 21 -- 26.

\z  E. Kolchin, Differential Algebra and Algebraic Groups, Academic Press, New
  York,  London, 1973.

\z  D. Kozlov, Combinatorial Algebraic Topology, Springer--Verlag Berlin, 2008.

\z  J. Ritt, Differential Algebra, Amer. Math. Soc. Colloq. Publ. 33, New York, 1950.

\z R. Stanley, Enumerative combinatorics, Vol. 2, Cambr, Univ. Press, 1999.

\z I. Trendafilov, D. Vladeva, The endomorphism semiring of a finite chain, Proc.  Techn. Univ.-Sofia, 61, 1, 2011, 9 -- 18.

 \z I. Trendafilov, D. Vladeva, Subsemirings of the endomorphism semi\-ring of a finite Chain, Proc.  Techn. Univ. - Sofia, vol. 61,  1, 2011, 19--28.

\z I. Trendafilov, Derivations in Some Finite Endomorphism Semirings, Discuss.   Math.
Gen. Algebra and Appl., Vol. 32, 2012, 77--100.

\z I. Trendafilov, Simplices in the Endomorphism Semiring of a Finite Chain,  Hindawi Publishing Corporation,
Algebra, Volume 2014, Art. ID 263605, 2014.

\z I. Trendafilov, D. Vladeva, Derivations in a triangle -- I and II part, Proc. Techn. Univ.-Sofia, 65, 1, 2015, 243 -- 262.

  \z D. Vladeva, I. Trendafilov,  Derivations in a tetrahedron, Appl. Math. in Eng. and Econ. -- 39th. Int. Conf. 2015,  American Institute of Physics Conference Proceedings 1690, 060004 (2015).

  \z D. Vladeva, I. Trendafilov,  Derivations in a tetrahedron -- I and II part, Proc. Techn. Univ.-Sofia, 65, 2, 2015, 117 -- 136.

 \z D. Vladeva, I. Trendafilov, A projection onto triangle of the 4--simplex is a derivation -- I and II part, Proc. Techn. Univ.-Sofia, 66, 3, 2016, 9 -- 28.

  \z D. Vladeva, I. Trendafilov, The Projections on Some Triangles of a Simplex  are Derivations,  Proc. of the Forty Five Spring Conference of the Union of Bulgarian Mathematicians, Pleven, April 6--10, 2016, 101 -- 107.

\z D. Vladeva,  Derivations in a endomorphism semiring, Serdica Math. J. 42, 2016, 251 -- 260.

\z J. Zumbr\"{a}gel, Classification of finite congruence-simple semirings with zero, J. Algebra Appl. 7, 2008,  363 -- 377.

\vspace{3mm}

{\bf \hspace{1mm} \textbf{Author:}} \hspace{1mm} Dimitrinka Vladeva, assoc. prof., Department "Mathematics and physics"$\hphantom{}$, LTU, Sofia, \emph{e-mail:}
d$\_$vladeva@abv.bg

\end{document}